\documentclass[11pt]{amsart}
\usepackage{amscd, amssymb}
\begin{document}

\newtheorem{thm}{Theorem}
\newtheorem{cor}[thm]{Corollary}
\theoremstyle{definition}
\newtheorem{definition}[thm]{Definition}
\newtheorem{example}[thm]{Example}

\numberwithin{equation}{section}

\newcommand{\supp}{\operatorname{supp}}
\newcommand{\Ind}{\operatorname{Ind}}
\newcommand{\tr}{\operatorname{tr}}

\title[Fell transformation groups]
{The transformation groups whose $C^*$-algebras\\ are Fell algebras}

\author[an Huef]{Astrid an Huef}

\address{Department of Mathematics\\
Dartmouth College\\
Hanover, NH 03755\\
USA}
\email{astrid.anhuef@dartmouth.edu}
\address{Department of Mathematics and Computer Science\\
University of Denver\\
Denver, CO 80208\\
USA
}
\email{astrid@cs.du.edu}
\thanks{This paper is based on part of the author's doctoral
  dissertation written at Dartmouth College. We thank Dana Williams
  for his supervision and help and Marc Rieffel for a 
  suggestion regarding the exposition of this paper.}

\date{July 23, 1999}

\begin{abstract}
Let $(G,X)$ be a locally compact  
transformation group, in which
$G$ acts freely on $X$. We show that the associated 
transformation-group $C^*$-algebra
$C_0(X)\rtimes G$ is  a
Fell algebra if  and only if $X$ is a Cartan $G$-space.
\end{abstract}

\maketitle

\def\N{{\mathbb N}}
\def\circle{{\mathbb T}}
\def\Z{{\mathbb Z}}
\def\C{{\mathbb C}}
\def\R{{\mathbb R}}
\def\P{{\mathbb P}}

\def\dach{^\wedge} 
\def\M{{\mathcal M}}
\def\calS{{\mathcal S}}

An irreducible representation $\pi$ of a $C^*$-algebra $A$ is a 
\emph{Fell point} of the spectrum $\hat A$ if there exist a
neighbourhood $V$
of $\pi$ in $\hat A$ and an element $a$ of the positive cone $A^+$ 
of $A$ such that $\sigma(a)$ is a
rank-one projection for all $\sigma\in V$. If every irreducible
representation is a Fell point then $A$ is a
\emph{Fell algebra}. 
A  $C^*$-algebra has \emph{continuous trace} if and only if its spectrum
is Hausdorff and it is a Fell algebra.

Let $(G,X)$ be a locally compact Hausdorff
transformation group: thus $G$ is a locally compact Hausdorff 
group and $X$ is a locally compact Hausdorff space together with a
jointly continuous map $(s,x)\mapsto
s\cdot x$ from $G\times X$ to $X$ such that $s\cdot (t\cdot x)=st\cdot
x$ and $e\cdot x=x$.
We assume that both  $G$ and $X$ are second countable.
The action of $G$ on $X$ lifts to a strongly continuous action
$\alpha$ of $G$ by automorphisms of $C_0(X)$ given by
$\alpha_s(f)(x)=f(s^{-1}\cdot x)$.
The associated transformation-group $C^*$-algebra $C_0(X)\rtimes G$ is
the $C^*$-algebra which is universal for the covariant representations
of the $C^*$-dynamical system $(C_0(X), G, \alpha)$ in the sense of
\cite{raeburn88:_takai}. More concretely, $C_0(X)\rtimes G$ is the
enveloping $C^*$-algebra of the Banach $*$-algebra $L^1(G, C_0(X))$ of 
functions $f:G\to C_0(X)$ which are integrable with respect to a
fixed left Haar measure $\mu$ on $G$
\cite[Section 7.6]{pedersen79:_automorphisms}, with
multiplication and involution  given by
$$
f*g(s):=\int_G f(r)\alpha_r(g(r^{-1}s))\ d\mu(r)\quad\text{and} 
\quad f^*(s):=\Delta(s^{-1})\alpha_s(f(s^{-1}))^*,
$$
where $\Delta$ is the modular function associated with $\mu$.

In \cite{green77:_orbit} Green showed that if the action of $G$ on $X$
is free, then   $C_0(X)\rtimes G$ has continuous
trace if and only if the action is proper.
The purpose of this note is to
obtain a parallel characterisation  of the  transformation 
groups whose $C^*$-algebras are Fell algebras.

\begin{thm}\label{sec:thm-main}
Let $(G, X)$ be a   second countable transformation 
group with $G$ acting freely on $X$. Then  $A=C_0(X)\rtimes G$ is a 
Fell algebra if and only if
$X$ is a Cartan $G$-space.
\end{thm}

A  subset $N$ of $X$ is {\emph{wandering}} if
the set $\lbrace s\in G: s\cdot N\cap N\neq\emptyset\rbrace$ is
relatively compact in $G$. 
If every point of $X$ has a wandering neighbourhood then $X$ is a
\emph{Cartan $G$-space} \cite[Definition 1.1.2]{palais61:_lie}.
When the action is free, these spaces are precisely the
principal bundles of Cartan by \cite[Theorem 1.1.3]{palais61:_lie},
and they are called principal $G$-bundles by Husemoller in 
\cite{husemoller75:_fibre}. 
The action of $G$ on $X$ is \emph{proper} (in the sense that the map 
$(s, x)\mapsto(s\cdot x, x)$ from $G\times X$ into $X\times X$ is a proper map)
if and only if every compact subset of $X$ is wandering.
If the action of $G$ on $X$ is proper then $X$ is a Cartan $G$-space,
but the converse is not true even if the action is free (see
Example~\ref{sec:palaisaction}).

For the proof of Theorem~\ref{sec:thm-main} we require some background
material.
In \cite{archbold93:_trans}, Archbold and Somerset
characterised Fell algebras $A$ as follows: define $S$ to be the set of
all $a\in A^+$  such that
$$
\tr(\pi_\alpha(a))\to\sum_{\pi\in L}\tr(\pi(a))<\infty\eqno (\dagger)
$$
whenever $\lbrace \pi_\alpha\rbrace$ is a properly convergent net in
$\hat A$
with limit set $L$. (A properly convergent net is a  convergent net 
which converges to every cluster point.) 
The linear span $\calS(A)$ of $S$ is a two-sided
(hereditary) ideal of $A$,  and  $A$ is a Fell algebra if and only if
$\calS(A)$ is  dense in  $A$ \cite[Corollary 3.7]{archbold93:_trans}.
To prove  Theorem~\ref{sec:thm-main} we show that if $X$ is not a
Cartan $G$-space then there is an element in the Pedersen ideal of an
ideal $I$ of $C_0(X)\rtimes G$ which is not an element of
$\calS(I)$. This suffices because an ideal of a Fell algebra is itself Fell,
and because the Pedersen ideal is minimal among dense
ideals \cite[Theorem 5.5.1]{pedersen79:_automorphisms}.

We thank Robert Archbold for pointing out that we can use
\cite[Theorem 3]{milicic73:_trace} instead of \cite[Theorem 2.6 and
Proposition 2.1 i)]{archbold97:_bded_tr_ideals} in the proof below.

\begin{proof}[Proof of Theorem~\ref{sec:thm-main}] 
We prove the easier direction first. 
If $X$ is a Cartan $G$-space then  the  orbits are  
closed by  \cite[Proposition 1.1.4]{palais61:_lie} and thus $A$ is CCR
\cite{williams81:_ccr}. Because the action is free the map 
$x\mapsto \pi_x:=\Ind_{\lbrace e\rbrace}^G(\epsilon_x)$ induces a
homeomorphism of $X/G$ onto $\hat A$  
\cite[Lemma 16]{green77:_orbit}. Fix $x\in X$ and let $U$ be an 
open wandering neighbourhood of $x$ in $X$.
The  action of $G$ on
$G\cdot U$ is proper by \cite[Proposition 1.2.4]{palais61:_lie}. 
Because the action is free
the algebra $C_0(G\cdot U)\rtimes G$ has continuous trace
\cite[Theorem 17]{green77:_orbit}, and since $C_0(G\cdot U)$ is a
$G$-invariant ideal of $C_0(X)$ this crossed product embeds naturally
as an ideal $J$ in $C_0(X)\rtimes G$ \cite[Lemma 16]{green77:_orbit}.
Viewing
$\pi_x$ as a representation of $J$ shows that it is a Fell point
of the open subset $\hat J$ of $\hat A$. Hence  $\pi_x$ is a 
Fell point of $\hat A$, and
$A$ is a Fell algebra.

Conversely, we show that if $X$ is not a Cartan $G$-space
then $A$ is not a Fell algebra. Some of the ideas 
come from  the proof of  \cite[Proposition 4.2]{williams81:_tr}.

Suppose there is a point $z$ in $X$ with no wandering
neighbourhood. To see that $A$ is not a Fell algebra it is enough to
produce an ideal $I$ of the form $C_0(Y)\rtimes G$ which is not a Fell
algebra. We may assume that $A$ has bounded trace, because otherwise
it is certainly not a Fell algebra.

To prove $I$ is not a Fell algebra, we construct a net of irreducible
representations $\lbrace\pi_{x_W}\rbrace$ of $\hat I$ parametrised by
a directed
family of neighbourhoods $W$ of $z$, which converges to its unique
cluster point $\pi_z$, and a positive element $d$ of the Pedersen ideal of $I$
such that $\tr(\pi_{x_W}(d))$ does not converge to $\tr(\pi_z(d))$.
The minimality of the Pedersen ideal among dense ideals of $I$ then
implies that the ideal $\calS(I)$ described in Equation~($\dagger$)
cannot be dense in $I$. It follows from
\cite[Corollary 3.7]{archbold93:_trans} that $I$ is not a Fell
algebra.
 
The construction of the element $d$ and the net
$\lbrace\pi_{x_W}\rbrace$ 
have to be
carefully intertwined to ensure that we can compare the traces by
finding appropriate eigenfunctions. We begin by building some
functions in $C_c(G\times X)$ which will later be used to define $d$.

Let $N$ be a compact
neighbourhood of $z$ in $X$ and  
$V_0$ an open symmetric neighbourhood of $e$ in $G$ with compact
closure $V$. Choose $f\in C_c(X)^+$ which is identically one on
$V_0\cdot N$. For $x\in X$ we define $f_x(s)=f(s\cdot x)$. Let
$F$ be a symmetric compact subset of $G$ such that $F\supseteq
\supp(f_z)$. Choose a non-negative function $b$ in $C_c(G\times X)$ such that
$b(r,x)=1$ whenever $r\in VF^2V$ and $x\in \supp(f)$. We can round off $b$ so that   $b(r,x)=0$ whenever
$(r,x)\notin V^2F^2V^2\times M$, where $M$ is a compact set 
containing the support of
$f$. 
We would like to use the function $(t,x)\mapsto
f(x)f(t^{-1}\cdot x)$ to define the element $d$, but this function
will in general not be compactly supported;  thus we will later
multiply it by the function $b$.

Let $K=V^2F^2V^3F$. 
By \cite[p. 61]{williams81:_tr} there exists 
an open neighbourhood $U$ of
$z$, such that for every $x\in U$, the support of ${f_x}_{|K}$ is contained in
$K\setminus V_0F$.
To see this, let $\phi:G\times X\to X$ be the map $(s,x)\mapsto s\cdot x$. 
Observe that $(K\setminus V_0F)\times \lbrace z\rbrace$ is compact and 
is contained in the open subset
$\phi^{-1}(X\setminus \supp(f))$ of
$G\times X$. 
For each $(s,z)\in (K\setminus V_0F)\times \lbrace z\rbrace$ 
there exist open neighbourhoods $V_s$ of $s$ and
$U_s$ of
$z$ such that $V_s\times U_s\subseteq \phi^{-1}(X\setminus
\supp(f))$. 
Choose a finite subcover $\lbrace
V_{s_i}\times U_{s_i}\rbrace_{i=1}^n$ and set
$U =\cap_{i=1}^n U_{s_i}$.
Then 
$$(K\setminus V_0F)\times\lbrace z\rbrace\subseteq
\cup_{i=1}^n V_{s_i} \times U \subseteq \phi^{-1}(X\setminus M).$$
Thus for every  $x\in U$, the function
$$
f^1_x(s) :=
\left\{\begin{array}{ll} f(s\cdot x)& \text{if}\,\,s\in K;
\\ 0 &\text{if}\,\, s\notin V_0F \end{array} \right.
$$
is  well-defined and continuous on $G$. This function will be an
eigenvector for the element $d$ to be constructed. 
Without loss of generality we may suppose $U\subseteq N$.

Since no neighbourhood of $z$ is wandering, for every
neighbourhood $W$ of $z$ there are $x_W\in W$ and
$s_W\notin K^{-1}K$ such that
 $s_W\cdot x_W\in W$.
Thus we get a net $\lbrace x_W\rbrace$ directed by the neighbourhoods
of $z$ which converges to $z$, 
and by passing to a subnet if necessary, we may
assume that this net is universal \cite[1.3.7]{pedersen89:_now} 
and that   $x_W\in U$ for all $W$.

Because $A$ has bounded trace and $\lbrace \pi_{x_W}\rbrace$ is a
convergent net, the limit set of
$\lbrace \pi_{x_W}\rbrace$ 
is discrete in the relative topology \cite[Theorem 3]{milicic73:_trace}. 
The map 
$x\mapsto \pi_x:=\Ind_{\lbrace e\rbrace}^G(\epsilon_x)$ induces a
homeomorphism of $X/G$ onto $\hat A$,
so there exists a neighbourhood $V$ of $G\cdot z$ in $X/G$ 
such that $G\cdot z$ is the unique limit point of
the net $\lbrace G\cdot x_W\rbrace$ in $V$.
We write $q:X\to X/G$ for the quotient map and let 
$Y= G\cdot (q^{-1}(V)\cap U)$. 
Since $Y$ is open we may assume that $\lbrace x_W\rbrace\subseteq Y$. Observe
that $\lbrace \pi_{x_W}\rbrace$ is the continuous image of a universal
net, and hence is universal; thus $\lbrace \pi_{x_W}\rbrace$ 
is properly convergent. 
 
We will show that the ideal
$I=C_0(Y)\rtimes G$ of $A$ is not a Fell algebra. 
Recall that our strategy to do this is to  construct
an element $d$ in the Pedersen ideal of
$I$ such that
$\lbrace \tr(\pi_{x_W}(d))\rbrace$ does not converge to  
$\tr(\pi_z(d))$. 
Let
$$B(t,x)=f(x)f(t^{-1}\cdot x)b(t^{-1},x)\Delta(t^{-1})^\frac{1}{2}
\quad\text{and}\quad D=\frac{1}{2}(B+B^*);$$
we will obtain the element $d$ by pushing enough of $D^**D$ into the
Pedersen ideal of $I$.

The   right Haar
measure $\nu$ on $G$ given by $\nu(E)=\mu(E^{-1})$ satisfies the hypothesis of
\cite[Lemma 4.14]{williams81:_ccr}, so $\pi_{x}$ is 
unitarily equivalent to the representation
$M_x\rtimes V$ on $L^2(G,\nu)$ given by
$$V(t)\xi(s)=\Delta(t)^\frac{1}{2}\xi(t^{-1}s)\quad\text{and}\quad
M_x(g)\xi(s)=g(s\cdot x)\xi(s)$$
 for $\xi\in L^2(G,\nu)$ and $g\in C_0(X)$.
Note that 
\begin{align*}
M_x\rtimes V(B)\xi(s) 
&=\int_G B(t,s\cdot x)\Delta(t)^\frac{1}{2}\xi(t^{-1}s) \ d\mu(t)\\
&=\int_G f(s\cdot x) f(u\cdot x) b(us^{-1}, s\cdot x) \xi(u) \ d\nu(u), 
\end{align*}
so that $M_x\rtimes V(B)$ is a Hilbert-Schmidt operator with kernel
$(s,u)\mapsto f(s\cdot x) f(u\cdot x) b(us^{-1},s\cdot x)$
in $C_c(G\times G)$.

By composing translation on the right by $(s_W)^{-1}$ with 
$f^1_{x_W}$ we see that 
$$
f^2_{x_W}(s) =
\left\{\begin{array}{ll} f(s\cdot x_W)& \text{if}\,\,s\in Ks_W;
\\ 0 &\text{if}\,\, s\notin V_0Fs_W \end{array} \right.
$$
is a well-defined continuous function on $G$ whenever $x_W\in U\cap Y$.
Note that $f_{x_W}^1$ and $f_{x_W}^2$ have disjoint supports because 
$s_W\notin K^{-1}K$.
The  $f^i_{x_W}$'s $(i=1,2)$ are  linearly
independent  eigenvectors for $M_{x_W}\rtimes V(D)$; 
as in the proof of \cite[Lemma 4.4]{williams81:_tr} one shows that
$$f(s\cdot x_W) f(u\cdot x_W)[b(us^{-1}, s\cdot x_W) + b(su^{-1}, u\cdot x_W)]
f_{x_W}^i(u)=2f^i_{x_W}(s) f(u\cdot x_W) f^i_{x_W}(u)$$ for all $s,u\in G$.
The eigenvalue corresponding to  $f^i_{x_W}$ is 
\begin{align*}
\lambda_{x_W}^i
 =\int_G f(u\cdot x_W) f_{x_W}^i(u) \ d\nu(u) 
=\int_G  f_{x_W}^i(u)^2 \ d\nu(u)
=\int_G  f_{x_W}^1(u)^2 \ d\nu(u)
=\lambda_{x_W}^1.
\end{align*}
If $x_W\in Y\cap U$ then $\lambda_{x_W}^1\geq\nu(V_0)>0$
because $u\in V_0$ implies $u\cdot x_W\in  V_0\cdot U\subseteq
V_0\cdot N$ 
and $f$ was chosen to be identically one on $V_0\cdot N$.
It follows that if $x_W\in Y\cap U$ then $M_{x_W}\rtimes V(D^**D)$  
is a positive
compact operator with  eigenvalue  $(\lambda^1_{x_W})^2$ greater or
equal to   $a:=\nu(V_0)^2$ of multiplicity at least $2$. 

On the other hand,  $M_z\rtimes V(D)$ is the rank-one operator
$f_z\otimes\bar f_z$ with eigenvalue
$\lambda_z^1=\|f_z\|^2=\int_G f_z(s)\  d\nu(s)=\lim\lambda_{x_W}^1>0$. 
Thus  $M_z\rtimes V(D^**D)$  is a rank-one operator with eigenvalue
$(\lambda_z^1)^2$.

Finally, we push enough of the operator $D^**D$ into the Pedersen
ideal of $I$ using the description in 
\cite[Theorem 5.5.1]{pedersen79:_automorphisms}.
Let $r$ be any non-negative  
function in $C_c((0,\infty))$ such that
$$
r(t) =
\left\{\begin{array}{ll} 0& \text{if}\,\,t<\frac{a}{3};
\\ 2t-\frac{2a}{3} &\text{if}\,\,\frac{a}{3}\leq t\leq\frac{2a}{3};
\\ t &\text{if}\,\, \frac{2a}{3}\leq t\leq \|D^**D\|. \end{array}
\right.
$$
Now $d=r(D^**D)$ is a positive element of the Pedersen ideal of $I$.  
If $x_W\in Y\cap U$ then
$\tr(M_{x_W}\rtimes V(d))\geq 2(\lambda_{x_W}^1)^2\to 2(\lambda_z^1)^2$ 
and $\tr(M_z\rtimes V(d))=(\lambda_z^1)^2$. The existence of such an
element implies that
$I$ and hence $A$ are not Fell algebras.
\end{proof}

A $C^*$-algebra $A$ always has a largest Fell ideal because the set of
Fell points is open in $\hat A$.

\begin{cor} Let $(G,X)$ be a  free second countable  transformation 
group such that 
$C_0(X)\rtimes G$ is Type {\rm I}. 
The largest Fell ideal of
$C_0(X)\rtimes G$ corresponds   to the open $G$-invariant subset 
$Y=\lbrace x\in X:\exists \text{\ a\  wandering\ neighbourhood\ of\ } x\rbrace$
of $X$.
\end{cor}

\begin{proof}
If the action is free and $C_0(X)\rtimes G$ is Type {\rm I} then every
closed ideal is of the form $C_0(Y)\rtimes G$ for some open
$G$-invariant subset of $X$ by \cite[Lemma 16]{green77:_orbit}. Now apply
Theorem~\ref{sec:thm-main} to $C_0(Y)\rtimes G$.
\end{proof}

The following transformation group is interesting for several
reasons. Its $C^*$-algebra is a Fell algebra which does not have
continuous trace. In addition, it has many ideals with continuous trace
but no largest continuous-trace ideal.
Recall that 
a positive element $a$ of a $C^*$-algebra $A$  is a 
\emph{continuous-trace element} if the
function $\pi\mapsto\tr(\pi(a))$ is finite and continuous on $\hat A$;
the $C^*$-algebra has continuous trace if the linear span
$\M(A)$ of these elements is a dense (hereditary) ideal of $A$.
In particular, the closure of $\M(C_0(X)\rtimes G)$ in the example
below is not a maximal continuous-trace ideal. 

\begin{example}\label{sec:palaisaction}
This transformation group was described in \cite[p. 298]{palais61:_lie}.
Let $X$ be the strip 
$\lbrace(x,y):-1\leq x\leq 1, y\in\R\rbrace$ in the plane. 
The group action is by $G=\R$.
Beyond the strip $-1<x<1$ the action moves a point according to
$$
t\cdot(1,y)=(1,y+t)\quad\text{and}\quad t\cdot(-1,y)=(-1, y-t).
$$
If $|x_0|<1$ let $C_{(x_0,y_0)}$ be the vertical 
translate of the graph of the equation $y=\frac{x^2}{1-x^2}$
which passes through $(x_0,y_0)$. 
Define $t\cdot(x_0,y_0)$ to be the point $(x,y)$ on $C_{(x_0,y_0)}$ 
such that the length of the arc of $C_{(x_0,y_0)}$ between 
$(x_0,y_0)$ and $(x,y)$ is $|t|$, and $x-x_0$ has the same sign as $t$.
That is, $(x_0,y_0)$ moves counter-clockwise along $C_{(x_0,y_0)}$ 
at unit speed.

A compact neighbourhood $N$ in $X$ is 
wandering if and only if it meets at most one boundary
line $x=1$ or $x=-1$;
thus $X$ is a Cartan $G$-space, but the action of 
$G$ on $X$ is not proper \cite[p. 298]{palais61:_lie}.
Since the action is free, $C_0(X)\rtimes G$ does not have continuous
trace, but it is a Fell algebra by  Theorem~\ref{sec:thm-main}.

The action of $G$ on both  
$Y=\lbrace (x,y): -1<x\leq 1, y\in\R\rbrace$ and 
$Z=\lbrace(x,y): -1\leq x<1, y\in\R\rbrace$
is free and proper, so that 
$I_1=C_0(Y)\rtimes G$ and  $I_2=C_0(Z)\rtimes G$
are distinct maximal ideals which both have continuous trace. 

We now identify the closure of the ideal of the continuous-trace
elements.
Let $U$ be any compact wandering neighbourhood of $(1,0)$ in $X$. 
Note that $U$ wandering means that $U$ does not meet
the line $x=-1$. The saturation $G\cdot U$ is not 
closed because we  can easily choose a sequence in $G\cdot U$
converging to $(-1,0)$.  The same is true for any compact 
wandering neighbourhood of $(1,y)$ or $(-1,y)$.
By \cite[Corollary 18]{green77:_orbit} the closure of the 
continuous-trace ideal $\M(C_0(X)\rtimes G)$
is $C_0(W)\rtimes G$ where $W$ consists of points $(x,y)$ which have a
wandering neighbourhood $U$ such that $G\cdot U$ is closed in
$X$. Thus
$W=\lbrace (x,y): -1<x<1,y\in\R\rbrace$
and  $\overline{\M(C_0(X)\rtimes G)}=C_0(W)\rtimes G=I_1\cap I_2$.
\end{example}

In \cite{aah:thesis} we characterise transformation groups whose
$C^*$-algebras have 
bounded trace. 

\bibliographystyle{amsplain}

 \end{document}